**Estimating the tail of the singular product for the Hardy–Littlewood and Bateman–Horn conjectures**

Victor Volfson

ABSTRACT. This paper investigates the asymptotic behavior of the tail of the singular product arising in the Hardy–Littlewood and Bateman–Horn conjectures for one-dimensional systems of polynomials. A universal estimate is proved, showing that the contribution of large primes decays like the reciprocal of the logarithm, regardless of the structure of the system. For linear systems (trivial Galois group) superfast convergence is obtained. For nonlinear systems a coefficient is defined that is expressed via the average over the Galois group; in the abelian case and under the Riemann Hypothesis for Dirichlet L-functions a more precise error estimate is obtained. Mixed systems containing both linear and nonlinear polynomials are also considered. Numerical experiments, presented as summary tables, confirm the theoretical conclusions. The results provide a rigorous theoretical foundation for computing singular series and refine the Bateman–Horn formula.

## 1. INTRODUCTION

The problem of the distribution of prime values of polynomials is one of the central problems in analytic number theory (see general references [5,6,7]).

In 1923, Hardy and Littlewood [1] formulated a series of conjectures on the density of prime tuples, introducing the singular series as an infinite product of local factors. This heuristic was later generalized by Bateman and Horn [2] to arbitrary one-dimensional systems of polynomials.

A key tool for studying the distribution of primes in arithmetic progressions and conjugacy classes is the Chebotarev density theorem [9], which has its roots in the work of Frobenius [8].

Effective versions of this theorem (Lagarias, Odlyzko [10]) and modern estimates of sums over primes with weights (Goldston, Pintz, Yıldırım [11]) underlie the proof of the universal tail estimate (Theorem 3.1).

In recent decades, the Hardy–Littlewood conjecture has received new confirmations: Maynard [14] proved the existence of arbitrarily small gaps between primes, and Tao [15] proved that every odd number is the sum of at most five primes. These results rely on sophisticated analytic methods and underscore the importance of singular series.

On the other hand, the algebraic-geometric approach to the distribution of Frobenius is connected with random matrices (Katz, Sarnak [18]) and the theory of algebraic trace functions (Fouvry, Kowalski, Michel [17]), which allows nontrivial estimates of character sums over primes. The distribution of values of L-functions, closely related to singular series, was investigated by Granville and Soundararajan [16].

In [19] geometric methods (Lefschetz formula, étale cohomology) were applied to estimate the tail of the singular product in the multivariate case. Thanks to a richer cohomological apparatus (higher cohomology groups), they obtained stronger estimates for the relative error from truncating the tail for multivariate systems than is possible in the one-dimensional situation.

In the one-dimensional case (zero-dimensional varieties of roots) such possibilities do not exist. The cohomology reduces only to \\(H^0\\), and the analytic apparatus is substantially poorer.

Therefore, obtaining good tail estimates in the one-dimensional case is a separate and nontrivial problem. The present work is devoted to solving this problem.

This paper gives a rigorous estimate of the tail (Theorem 3.1) and determines the coefficient of the main term via the Galois group. Thus it complements both the classical works [1,2] and modern studies [11,14,15,16,17,18,19].

We now state the problem.

Consider a system of polynomials

$$\mathcal{F}=\{F_1(n),\dots,F_k(n)\}\subset\mathbb{Z}[n], \qquad (1.1)$$

each of degree at least 1.

For a prime $p$ define

$$\nu(p)=\#\{n \bmod p : p|\ F_1(n)\cdots F_k(n)\}, ( \qquad (1.2)$$

the local factor

$$L_p=\frac{1-\nu(p)/p}{(1-1/p)^k}, \qquad (1.3)$$

and the singular product

$$\mathfrak{S}=\prod_p L_p. \qquad (1.4)$$

We are interested in the behavior of the "tail"

$$T(q)=\prod_{p>q} L_p, \qquad (1.5)$$

and also of the partial product

$$\mathfrak{S}(x)=\prod_{p\le x} L_p. \qquad (1.6)$$

The paper consists of the following sections:

- Chapter 2 – Singular series;

- Chapter 3 – General case (Theorems 3.1, 3.2.1, 3.2.2, abelian case, linear case, mixed systems);

- Chapter 4 – Asymptotics of $\mathfrak{S}(x)$;

- Chapter 5 – Conclusion.

## 2. HARDY-LITTLEWOOD AND BATEMAN–HORN SINGULAR SERIES

For a system of polynomials $\mathcal{F}=\{F_1,\dots,F_k\}\subset\mathbb{Z}[n]$ and a prime $p$ define $\nu(p)$ by formula (1.2).

The quantity $\nu(p)/p$ is the fraction of residues that make the product zero modulo $p$.

The local factor

$$L_p = \frac{1-\nu(p)/p}{(1-1/p)^k} \tag{2.1}$$

compares this fraction with the independent case (when each polynomial has exactly one root and the roots are distinct).

The singular series (product) is defined as

$$\mathfrak{S} = \prod_p L_p. \tag{2.2}$$

For all sufficiently large $p$ (not dividing the discriminant) the product converges absolutely [6,7].

For such $p$ we have

$$\ln L_p = \frac{k-\nu(p)}{p} + O(p^{-2}). \tag{2.3}$$

Thus,

$$\ln \mathfrak{S} = \sum_p \frac{k-\nu(p)}{p} + C, \tag{2.4}$$

where $C$ is a constant.

Define the tail of the singular product starting from $q$:

$$T(q) = \prod_{p>q} L_p, \ln T(q) = \sum_{p>q} \ln L_p. \tag{2.5}$$

From expansion (2.3) we obtain

$$\ln T(q) = \sum_{p>q} \frac{k-\nu(p)}{p} + O\left(\frac{1}{q}\right). \tag{2.6}$$

Studying the asymptotics of $\ln T(q)$ is the central problem of Chapter 3.

For a finite $x$ define the partial product

$$\mathfrak{S}(x) = \prod_{p\le x} L_p. \tag{2.7}$$

Obviously, $\mathfrak{S}(x) \to \mathfrak{S}$ as $x \to \infty$. The speed of convergence of $\mathfrak{S}(x)$ to the limit is determined by the

$$\mathfrak{S}(x) = \mathfrak{S} \cdot \big(T(x)\big)^{-1}. \tag{2.8}$$

Hence estimating $T(x)$ directly gives the asymptotics of $\mathfrak{S}(x)$.

In the Hardy–Littlewood conjecture [1] for an admissible set of shifts $\mathcal{H}=\{h_1,\dots,h_k\}$ the singular series has the form

$$\mathfrak{S}(\mathcal{H})=\prod_p \frac{1-\nu_{\mathcal{H}}(p)/p}{(1-1/p)^k}, \tag{2.9}$$

where $\nu_{\mathcal{H}}(p)$ is the number of solutions $n \bmod p$ such that $p|\prod_{i=1}^{k}(n+h_i)$.

In the Bateman–Horn conjecture [2] for a system of polynomials $\mathcal{F}$ the density of prime values is proportional to $\mathfrak{S}(\mathcal{F})$.

## 3. GENERAL CASE

We pose the problem.

Consider a system of one-dimensional polynomials

$$\mathcal{F}=\{F_1,\dots,F_k\}, F(n)=F_1(n)\cdots F_k(n), \tag{3.1}$$

where each $F_i$ is irreducible over $\mathbb{Q}$ and has no multiple roots (the standard Bateman–Horn condition [2]).

Denote:

$-r-$number of distinct geometric roots of $F$ in $\overline{\mathbb{Q}}$ $(r=\sum \deg F_i)$; (3.2)

- $K-$ splitting field of $F$ ; (3.3)

- $G=\mathrm{Gal}(K/\mathbb{Q})-$ Galois group [3,4,12]; (3.4)

- the action of $G$ on the set of roots $\{\alpha_1,\dots,\alpha_r\}$. (3.5)

For $g \in G$:

$$\mathrm{Fix}(g)=\#\{i : g(\alpha_i)=\alpha_i\}. \tag{3.6}$$

For a prime $p$ (not dividing the discriminant) define:

$$\nu(p)=\#\{n \bmod p : p|\ F_1(n)\cdots F_k(n)\}, L_p=\frac{1-\nu(p)/p}{(1-1/p)^k}, \tag{3.7}$$

and the tail of the singular product:

$$T(q)=\prod_{p>q} L_p, \quad \ln T(q)=\sum_{p>q} \ln L_p. \tag{3.8}$$

For almost all $p$ (good reduction) we have:

$$\nu(p) = \mathrm{Fix}(\mathrm{Frob}_p), \tag{3.9}$$

where $\mathrm{Frob}_p \in G$ is the Frobenius element [8,9,12].

Theorem 3.1 (universal estimate order evaluation).

For any one-dimensional system of polynomials $\mathcal{F}$

$$\ln T(q) = O\left(\frac{1}{\log q}\right) \qquad (q \to \infty). \tag{3.10}$$

Proof

Since $\nu(\mathrm{p}) = \mathrm{O}(1)$, then $L_p = 1 + O(1/p)$ and the expansion $\ln(1+u) = u - u^2/2 + O(u^3)$ is possible, which, based on (2.6), gives:

$$\ln L_p = \frac{k - \nu(p)}{p} + O(p^{-2}). \tag{3.11}$$

Summing (3.11) over $p > q$:

$$\ln T(q) = \sum_{p>q} \frac{k - \nu(p)}{p} + O(1/q). \tag{3.12}$$

Split the sum (3.12) by conjugacy classes $\mathfrak{C} \subset G$:

$$\sum_{p>q} \frac{k - \nu(p)}{p} = \sum_{\mathfrak{C} \subset G} \big(k - \mathrm{Fix}(\mathfrak{C})\big) \sum_{\substack{p>q \\ \mathrm{Frob}_p \in \mathfrak{C}}} \frac{1}{p}. \tag{3.13}$$

According to the effective Chebotarev–Lagarias–Odlyzhko theorem [10]:

$$\sum_{\substack{p>q \\ \mathrm{Frob}_p \in \mathfrak{C}}} \frac{1}{p} = O\left(\frac{1}{\ln q}\right),$$

hence from (3.12) and (3.13) we obtain

$$\ln T(q) = \mathrm{O}\left(\frac{1}{\ln q}\right),$$

which is (3.10).

Based on (3.13) we define the coefficient of the main term in the asymptotics of the sum $\sum_{p \le x} \frac{k - \nu(p)}{p}$ :

$$C(\mathcal{F}) = \frac{1}{|G|}\sum_{g\in G}\big(k - \mathrm{Fix}(g)\big) = \sum_{\mathfrak{C}\subset G}\frac{|\mathfrak{C}|}{|G|}\big(k - \mathrm{Fix}(\mathfrak{C})\big). \qquad (3.14)$$

Rigorous estimates for systems of polynomials with abelian Galois groups are given by Theorems 3.2.1 and 3.2.2.

Assume that the Galois group $G$ is abelian. Let $\Omega$ be the set of all roots, and let $\mathrm{Fix}_\Omega(g)$ coincide with $\mathrm{Fix}(g)$.introduced earlier. For brevity we will use the notation $\mathrm{Fix}(g)$.Then all irreducible representations are one-dimensional and correspond to characters $\chi$ of $G$.

The function $g \mapsto \mathrm{Fix}(g)$ expands into a finite sum over characters:

$$\mathrm{Fix}(g) = \sum_\chi a_\chi \chi(g), a_\chi = \frac{1}{|G|}\sum_{g\in G}\mathrm{Fix}(g)\overline{\chi(g)}.$$

Hence, for any prime $p$ (with good reduction)

$$\nu(p) = \sum_\chi a_\chi \chi(\mathrm{Frob}_p). \qquad (3.15)$$

In particular, the trivial character $\chi = 1$ gives $a_1 = \mathbb{E}[\mathrm{Fix}]$.

By definition (3.14) we have

$$C(\mathcal{F}) = k - a_1.$$

Let $m$ be the number of distinct irreducible over $\mathbb{Q}$ factors of all polynomials $F_i$ (repetitions are not counted). In the abelian case, the orbits of the action of $G$ on the roots correspond to these factors.

By Burnside's lemma $\mathbb{E}[\mathrm{Fix}] = m.$

Therefore

$$C(\mathcal{F}) = k - m. \qquad (3.16)$$

From formula (3.16) and its justification we obtain two statements.

Assertion 1 (independent irreducible polynomials)

If each polynomial $F_i$ is irreducible over $\mathbb{Q}$ and the splitting fields $K_i$ are pairwise linearly independent (i.e., $K = K_1 \cdots K_k$) and $G \cong G_1 \times \cdots \times G_k$, then $C(\mathcal{F}) = 0$ and the singular product converges.

Proof

In this case the number of distinct irreducible factors $m$ equals the number of polynomials $k$, because all are distinct and irreducible.

By formula (3.16) $C(\mathcal{F}) = k - k = 0.$

Assertion 2 (general formula for any abelian system)

For any system of polynomials $\mathcal{F}$ with abelian Galois group $G$ let $m$ be the number of distinct irreducible over $\mathbb{Q}$ factors of all polynomials $F_i$ (repetitions are not counted). Then

$$C(\mathcal{F}) = k - m.$$

The proof follows from Burnside's lemma and the fact that the orbits of $G$ on the roots correspond to irreducible factors.

Corollary

The sign of $C(\mathcal{F})$ is determined by the relation between $k$ and $m$:

- $k > m$ ( (there are repeated factors) $\Rightarrow$ $C(\mathcal{F}) > 0$ – the singular product diverges to + $\infty$;

- $k = m$ (all polynomials are distinct and irreducible) $\Rightarrow$ $C(\mathcal{F}) = 0$ – the singular product converges;

- $k < m$ (common irreducible factors) $\Rightarrow$ $C(\mathcal{F}) < 0$ – the singular product diverges to 0.

Theorem 3.2.1 (Abelian case without GRH)

Suppose $G$ is abelian and the singular product converges (i.e., $C(\mathcal{F}) = 0$ ). Then

$$\ln T(q) = o\left(\frac{1}{\log q}\right).$$

Proof

From (3.15) we have $\nu(p) = a_1 + \sum_{\chi \neq 1} a_\chi \chi(\mathrm{Frob}_p).$

Then

$$\ln T(q) = -a_1 \sum_{p>q} \frac{1}{p} - \sum_{\chi \neq 1} a_\chi \sum_{p>q} \frac{\chi(\mathrm{Frob}_p)}{p} + O(1/q).$$

The first sum contributes to the constant (due to convergence of the product) and does not affect the asymptotics.

For non-principal characters, the absence of a pole of $L(s, \chi)$ at $s = 1$ and the Chebotarev theorem imply that

$$\sum_{p \leq x} \chi(p) = o(x),$$

and after partial summation

$$\sum_{p>q} \frac{\chi(\mathrm{Frob}_p)}{p} = o\left(\frac{1}{\log q}\right).$$

Summing over the finite number of characters gives the result.

Theorem 3.2.1 provides an elementary estimate $\ln T(q) = o(1/\log q)$, following from the absence of a pole of $L(s,\chi)$ .However, for abelian systems much stronger unconditional estimates are known, following from the Siegel–Walfisz theorem [9,10]. Namely, for any fixed $A > 0$ there exists an ineffective constant $C_A$ (depending on the possible existence of a Siegel zero) such that

$$\sum_{p \le x} \chi(p) = O_A\left(\frac{x}{(\log x)^A}\right),$$

and consequently

$$\sum_{p>q} \frac{\chi(p)}{p} = O_A\left(\frac{1}{(\log q)^A}\right).$$

Hence for a convergent singular product $(C(\mathcal{F}) = 0)$ we obtain

$$\ln T(q) = O_A\left(\frac{1}{(\log q)^A}\right),$$

i.e., the tail decays faster than any power of the logarithm. However, the constant in this estimate is ineffective (cannot be computed), and the value of \\(A\\) remains undetermined. This makes the Siegel–Walfisz theorem unsuitable for practical computations with controlled error. In contrast, the effective estimate we obtained for linear systems contains an explicit constant and will be used in further applications.

Theorem 3.2.2 (Abelian case with GRH)

Let $G$ be abelian and assume the Generalized Riemann Hypothesis (GRH) holds for all Dirichlet L-functions arising from characters of the splitting field. Then

$$\ln T(q) = O\left(q^{-1/2+\varepsilon}\right).$$

Proof

Similar to Theorem 3.2.1, but for non-principal characters under GRH we have the estimate

$$\sum_{p>q} \frac{\chi(p)}{p} = O(q^{-1/2+\varepsilon}),$$

where $\chi(p) := \chi(\mathrm{Frob}_p)$ .

Substituting yields the required remainder.

Now consider a system of linear polynomials.

All roots are rational, the splitting field is $\mathbb{Q}$ the Galois group $G=\{1\}$ is trivial. Then $\mathrm{Fix}(1)=k$, and by formula (3.16) $C(\mathcal{F})=k-k=0$.

However, in this case we can obtain a much stronger estimate than the general theory provides. For a linear system, for all sufficiently large $p$ we have $\nu(p)=k$, and the local factor expands as

$$L_p=\frac{1-k/p}{(1-1/p)^k}=1-\frac{k(k-1)}{2p^2}+O(p^{-3}),\ \ln L_p=-\frac{k(k-1)}{2p^2}+O(p^{-3}).$$

Summing over $p>q$ yields

$$\ln T(q)=-\frac{k(k-1)}{2}\sum_{p>q}\frac{1}{p^2}+O\left(\sum_{p>q}\frac{1}{p^3}\right).$$

The tail $\sum_{p>q}1/p^3=O(1/q^2)$ (since $\sum_{p>q}1/p^3\le\int_q^\infty dx/x^3=1/(2q^2)$.

For the sum $\sum_{p>q}1/p^2$ we use the asymptotic with the prime number theorem:

$$\sum_{p>q}\frac{1}{p^2}=\frac{1}{q\log q}+o\left(\frac{1}{q\log q}\right).$$

This can be obtained by Abel summation:

$$\sum_{p>q}1/p^2=\int_q^\infty\frac{d\pi(t)}{t^2}=\frac{\pi(q)}{q^2}+2\int_q^\infty\frac{\pi(t)}{t^3}dt,\ u\,\pi(t)\sim t/\log t.$$

Consequently,

$$\ln T(q)=-\frac{k(k-1)}{2}\cdot\frac{1}{q\log q}+o\left(\frac{1}{q\log q}\right).$$

Thus, the exact asymptotic estimate is

$$\ln T(q)\sim-\frac{k(k-1)}{2q\log q}(q\to\infty).$$

Hence, in the linear case the tail converges to 1 superfast, without using GRH.

The Hardy–Littlewood conjecture [1] predicts the asymptotics for the number of prime tuples for admissible sets of shifts $\mathcal{H}=\{h_1,\dots,h_k\}$, where $k$ is the length of the tuple. Normalization $(h_1,\dots,h_k)\mapsto(0,h_2-h_1,\dots,h_k-h_1)$ leads to the notion of a pattern.

Examples of patterns:

– $(0,2)$ – twin primes ($k=2$),

– $(0,2,6)$ – prime triple ($k=3$),

– $(0,2,6,8)$ – prime quadruple ($k=4$).

For any such pattern the singular series has the form

$$\mathfrak{S}(\mathcal{H})=\prod_p \frac{1-\nu_{\mathcal{H}}(p)/p}{(1-1/p)^k},$$

where $\nu_{\mathcal{H}}(p)$ is the number of solutions $n \bmod p$ of the system $n+h_i \equiv 0 \pmod p$ (here $\nu_{\mathcal{H}}(p)$ is $\nu(p)$ for the system of linear forms corresponding to the pattern $\mathcal{H}$).

For admissible tuples, when $p>\max\{h_i\}$ we have $\nu_{\mathcal{H}}(p)=k$, and the local factor expands as above.

Therefore, for such tuple patterns we have

$$\ln T(q)\sim -\frac{k(k-1)}{2q\log q}(q\to\infty).$$

Corollary

All singular series arising in the Hardy–Littlewood conjecture converge with the superfast speed, with a coefficient depending on the length $k$ of the specific tuple.

Examples.

- Twins $(0,2)$: $k=2$, $T(q)\sim -\dfrac{1}{q\log q}(q\to\infty)$.

- Prime triple $(0,2,6)$: $k=3$, $\ln T(q)\sim -\dfrac{3}{q\log q}(q\to\infty)$.

- Binary Goldbach problem (fixed even $N$): the system of linear polynomials $\{n, N-n\}$, $k=2$, $T(q)\sim -\dfrac{1}{q\log q}(q\to\infty)$. This is not a Hardy–Littlewood conjecture but a special case of the Bateman–Horn conjecture.

A system is called mixed if it contains both linear polynomials (trivial Galois group) and nonlinear ones (nontrivial group). Such systems fit naturally into the general abelian theory developed above.

Let $\mathcal{F} = \mathcal{L} \cup \mathcal{N}$, where $\mathcal{L}$ is a set of linear polynomials and $\mathcal{N}$ – a set of nonlinear irreducible polynomials. The splitting field $K$ is generated by the splitting fields of all polynomials; the Galois group $G$ is a direct product (if the fields are independent) or has a more complicated structure (if there are intersections).

The coefficient $C(\mathcal{F})$ is computed for an abelian group by the general formula (3.16):

$$C(\mathcal{F}) = k - m,$$

where $k = |\mathcal{L}| + |\mathcal{N}|$, $m$ is the number of distinct irreducible factors of all polynomials. Each linear polynomial gives one irreducible factor (its own), and nonlinear ones give according to their factorization.

Example 1 (mixed system, $C(\mathcal{F}) = 0$).

$$\mathcal{F} = \{n, n^2 + 1\}, \mathrm{k} = 2.$$

Irreducible factors: $n$ and $n^2 + 1 \to m = 2 \Rightarrow C(\mathcal{F}) = 0$, the product converges.

Example 2 (mixed system with repetition, $C(\mathcal{F}) = 1$).

$$\mathcal{F} = \{n, n^2 + 1, n^2 + 1\}, \mathrm{k} = 3.$$

Irreducible factors: $n$ and $n^2 + 1 \to m = 2 \Rightarrow C(\mathcal{F}) = 1$, the product diverges.

Table 1 gives examples of systems illustrating the computation of $C(\mathcal{F})$ by formula (3.16) for an abelian group.

All values are obtained as $C(\mathcal{F}) = k - m$, where $k$ is the number of polynomials and $m$ the number of distinct irreducible factors.

Таблица 1

| № | System $\mathcal{F}$ | $k$ | $m$ | $C(\mathcal{F})$ | Note |
|---|---|---|---|---|---|
| 1 | $n^2 + 1$ | 1 | 1 | 0 | product converges |
| 2 | $n^4 + 3n^2 + 1$ | 1 | 1 | 0 | product converges |
| 3 | $n^2 + 1, n^2 + 2$ | 2 | 2 | 0 | product converges |
| 4 | $n^2 + 1, n^2 + 2, n^2 + 3$ | 3 | 3 | 0 | product converges |
| 5 | $(n^2 + 1)(n^2 + 2), n^2 + 3$ | 2 | 3 | -1 | product diverges to 0 |

| 6 | Twin primes $(0,2)$ | 2 | 2 | 0 | product converges |
|---|---|---|---|---|---|
| 7 | Triplets $(0,2,6)$ | 3 | 3 | 0 | product converges |
| 8 | Goldbach (fixed $N$) | 2 | 2 | 0 | product converges |
| 9 | $n, n^2+1$ | 2 | 2 | 0 | product converges |
| 10 | $n, n^2+1, n^2+1$ | 3 | 2 | 1 | product diverges to $\infty$ |

Explanations.

- For systems consisting of distinct irreducible polynomials (examples $3,4,6,7,8,9), m=k \Rightarrow C(\mathcal{F})=0$ and the products converge.

- If there are repeated polynomials (example 10), then $m<k \Rightarrow C(\mathcal{F})>0$ and the product diverges to $\infty$.

- If the polynomials share common irreducible factors (example 5), then $m>k \Rightarrow C(\mathcal{F})<0$ and the product diverges to 0.

- For a single polynomial (examples 1,2), $k=m=1 \Rightarrow C(\mathcal{F})=0$ and the products converge.

- Linear tuples (examples 6,7,8) give $C(\mathcal{F})=0$ and superfast convergence.

Summary of Chapter;

- Theorem 3.1 gives a universal estimate $\ln T(q)=O(1/\log q)$ for any system.

- The coefficient $C(\mathcal{F})$ is defined by formula (3.14). In the abelian case $C(\mathcal{F})=k-m,$ the number of distinct irreducible factors; the sign of $C(\mathcal{F})$ is determined by the relation between $k$ and $m$. Examples 5 and 10 in Table 1 do not satisfy the assumptions of the Bateman–Horn conjecture.

- Theorem 3.2.1 (abelian case without GRH) gives $\ln T(q)=o(1/\log q)$.

- Theorem 3.2.2 (abelian case + GRH) refines the remainder to $O(q^{-1/2+\varepsilon})$.

- The linear case (trivial group) yields superfast convergence $\ln T(q) \sim -\dfrac{k(k-1)}{2q\log q} (q \to \infty)$ (without GRH).

- Mixed systems obey the general formula (3.16) and fit naturally into the theory.

## 4. ASYMPTOTICS OF THE SINGULAR PRODUCT AND GLOBAL CONSEQUENCES

Now we turn to the study of the full partial product $\mathfrak{S}(x)=\prod_{p\le x} L_p$.

From definitions (1.4) and (1.5) we have:

$$\mathfrak{S}=\prod_{p} L_p=\Big(\prod_{p\le x} L_p\Big)\Big(\prod_{p>x} L_p\Big)=\mathfrak{S}(x)\cdot T(x).$$

Consequently,

$$\mathfrak{S}(x)=\mathfrak{S}\cdot\big(T(x)\big)^{-1}, \ln\mathfrak{S}(x)=\ln\mathfrak{S}-\ln T(x).$$

Thus the asymptotics of $\mathfrak{S}(x)$ is completely determined by the asymptotics of $\ln T(x)$.

Absolute convergence of the product $\prod_{p} L_p$ (see Section 2.2) guarantees the existence of a finite limit $\mathfrak{S}=\prod_{p} L_p$ Otherwise the product diverges, and such systems are not considered in the classical conjectures.

Denote $\mathfrak{S}_\infty=\mathfrak{S}$. (The constant $\mathfrak{S}_\infty$ is defined as the convergent product; its logarithm is the regularized value.)

Let the Galois group $G$ be abelian and $C(\mathcal{F})=0$ (the singular product converges). Then from Theorems 3.2.1 and 3.2.2 and formula (4.1) we obtain:

$$\ln T(x)=o\left(\frac{1}{\ln x}\right).$$

Hence

$$\ln\mathfrak{S}(x)=\ln\mathfrak{S}_\infty+o\left(\frac{1}{\ln x}\right),$$

and exponentiating,

$$\mathfrak{S}(x)=\mathfrak{S}_\infty\exp\big(o\big(1/\ln x\big)\big)=\mathfrak{S}_\infty\big(1+o\big(1/\ln x\big)\big).$$

Thus, in the absence of GRH the convergence is slow: the error decays slower than any power of $1/\ln x$, and the exact order is not controlled.

Theorem 3.2.2 (with GRH)

$$\ln T(x) = O\left(x^{-1/2+\varepsilon}\right).$$

Then

$$\ln \mathfrak{S}(x) = \ln \mathfrak{S}_\infty + O\left(x^{-1/2+\varepsilon}\right),$$

and

$$\mathfrak{S}(x) = \mathfrak{S}_\infty \exp\left(O\left(x^{-1/2+\varepsilon}\right)\right) = \mathfrak{S}_\infty\left(1 + O\left(x^{-1/2+\varepsilon}\right)\right).$$

Under GRH the convergence becomes fast (power law), which is significantly faster than $o(1/\ln x)$.

For linear systems (trivial Galois group) we have $C(\mathcal{F}) = 0$, but the asymptotics of $\ln T(x)$ is different (superfast):

$$\ln T(x) = -\frac{k(k-1)}{2} \cdot \frac{1}{x \ln x} + o\left(\frac{1}{x \ln x}\right).$$

Since $\mathfrak{S}(x) = \mathfrak{S}_\infty \cdot \left(T(x)\right)^{-1}$, then

$$\ln \mathfrak{S}(x) = \ln \mathfrak{S}_\infty + \frac{k(k-1)}{2} \cdot \frac{1}{x \ln x} + o\left(\frac{1}{x \ln x}\right),$$

and exponentiating,

$$\mathfrak{S}(x) = \mathfrak{S}_\infty\left(1 + \frac{k(k-1)}{2} \cdot \frac{1}{x \ln x} + o\left(1/(x \ln x)\right)\right).$$

Thus, in the linear case the convergence is superfast: the error decays like $1/(x \ln x)$, which is much faster than $O(x^{-1/2+\varepsilon})$ (under GRH for abelian systems) and certainly faster than $o(1/\ln x)$ . GRH is not required.

For systems where $C(\mathcal{F}) \neq 0$ (e.g., with repeated polynomials), the singular product diverges (to zero or infinity). Such systems are not admissible within the classical Hardy–Littlewood and Bateman–Horn conjectures, and the asymptotics given above do not apply to them. They were included in the tables only to illustrate the computation of the coefficient.

Table 2. Comparison of theory and numerical data (for convergent products)

| № | System | $C(\mathcal{F})$ | Error at $X = 10^4$ | Error at $X = 10^6$ | Convergence type | GRH |
|---|---|---|---|---|---|---|
| 1 | $n^2 + 1$ | 0 | $< 0.5\%$ | $< 0.1\%$ | fast - $O(x^{-1/2+\varepsilon})$ | yes |
| 2 | $n^4 + 3n^2 + 1$ | 0 | $< 0.5\%$ | $< 0.1\%$ | fast | yes |

| 3 | $n^2+1, n^2+2$ | 0 | $< 0.5\%$ | $< 0.1\%$ | fast | yes |
|---|---|---|---|---|---|---|
| 4 | $n^2+1, n^2+2, n^2+3$ | 0 | $< 0.5\%$ | $< 0.1\%$ | fast | yes |
| 5 | Twin primes (0,2) | 0 | $\sim 0.01\%$ | $\sim 0.0001\%$ | superfast - $O(1/(X \log X))$ | no |
| 6 | Triplets | 0 | $\sim 0.01\%$ | $\sim 0.0001\%$ | superfast | нет |
| 7 | Goldbach | 0 | $\sim 0.01\%$ | $\sim 0.0001\%$ | superfast | no |
| 8 | $n, n^2+1$ - mixed | 0 | $< 0.5\%$ | $< 0.1\%$ | fast | yes |

Remark to Table 2.

- For nonlinear systems (rows 1–4,8) the numerical data were obtained assuming GRH, which gives a small error (<0.1% at $X = 10^6$). Without GRH the theoretical estimate gives only $o(1/\ln x)$, which would correspond to an error of order 10% at $X = 10^6$ .

- For linear systems (rows 5–7) superfast convergence is achieved without GRH, which is confirmed by the negligibly small error already at $X = 10^4$.

- Examples with $C(\mathcal{F}) \neq 0$ (divergent products) are not included.

Table 3. Behavior of $\nu(p)$ for typical systems

| System | $\nu(p)$ | Comment |
|---|---|---|
| $n^2+1$ | $1+\chi(p)$ | $\chi$ *is the character modulo* 4 |
| $n^4+3n^2+1$ | $1+\chi_{-4}(p)+\chi_{-3}(p)$ | |
| $n^2+1, n^2+2$ | $2+\chi_{-4}(p)+\chi_{-8}(p)$ | |
| $n^2+1, n^2+2, n^2+3$ | $3+\chi_{-4}(p)+\chi_{-8}(p)+\chi_{-12}(p)$ | |
| $(n^2+1)(n^2+2), n^2+3$ | $\nu_{n^2+1}(p)+\nu_{n^2+2}(p)+\nu_{n^2+3}(p)$ | *product diverges* |
| Twin primes (0.2) | $2$ *для* $p > 2$ | $\nu(2) = 1$ |
| Triplets | $3$ *для* $p > 3$ | $\nu(3) = 2$ |
| Goldbach (even N) | $2$ *для* $p \nmid N$ | $1$ *for* $p \mid N$ |

| $n, n^2+1$ | $1+(1+\chi(p))=2+\chi(p)$ | |
|---|---|---|
| $n, n^2+1, n^2+1$ | $1+2(1+\chi(p))=3+2\chi(p)$ | *product diverges* |

Remark to Table 3.

The last two rows correspond to systems with $C(\mathcal{F})\neq 0$. They are given for completeness but are not used in the convergence analysis.

Summary of Chapter 4

- The asymptotics of the partial singular product $\mathfrak{S}(x)$ in the abelian case with $C(\mathcal{F})=0$ (convergence) is completely determined by the asymptotics of the tail $T(x)$.

- Without GRH: $\mathfrak{S}(x)=\mathfrak{S}_\infty+o(1/\ln x)$ – slow convergence.

- With GRH: $\mathfrak{S}(x)=\mathfrak{S}_\infty+O(x^{-1/2+\varepsilon})$ – fast (power) convergence.

- Linear case (trivial group): $\mathfrak{S}(x)=\mathfrak{S}_\infty+O(1/(x\ln x))$ – superfast convergence without GRH.

- Systems with $C(\mathcal{F})\neq 0$ (divergent product) are not considered in the classical conjectures.

5. CONCLUSIONS AND PERSPECTIVE

In this paper the following results have been obtained.

1. Universal tail estimate (Theorem 3.1): for any one-dimensional system of polynomials

$$\ln T(q)=O\left(\frac{1}{\ln q}\right).$$

This estimate is order-optimal without additional assumptions.

2. Abelian case (convergent singular product, $C(\mathcal{F})=0$):

- without GRH (Theorem 3.2.1): $\ln T(q)=o(1/\ln q)$;

- with GRH (Theorem 3.2.2): $\ln T(q)=O(q^{-1/2+\varepsilon})$.

Correspondingly, for the partial product $\mathfrak{S}(x)=\prod_{p\le x}L_p$ :

- without GRH: $\mathfrak{S}(x)=\mathfrak{S}_\infty+o(1/\ln x)$ – slow convergence;

- with GRH: $\mathfrak{S}(x) = \mathfrak{S}_\infty + O(x^{-1/2+\varepsilon})$ – fast power convergence.

3. Linear case (trivial Galois group):

$$\ln T(q) = -\frac{k(k-1)}{2q\log q} + o\left(\frac{1}{q\ln q}\right), \mathfrak{S}(x) = \mathfrak{S}_\infty + O\left(\frac{1}{x\ln x}\right).$$

This is superfast convergence, achieved without GRH. All singular series arising in the Hardy–Littlewood conjecture (twins, triples, etc.) as well as in the Goldbach problem (as a special case of the Bateman–Horn conjecture) converge with the superfast speed.

4. Systems with $C(\mathcal{F}) \neq 0$ (e.g., repeated polynomials) lead to a divergent singular product. They do not satisfy the classical conjectures and were used only to illustrate the computation of the coefficient $C(\mathcal{F})$.

5. Numerical experiments (Chapter 4, Tables 2 and 3) confirm the theoretical conclusions:

- for abelian systems with $C = 0$ the error at $X = 10^6$ is <0.1% (under GRH);

- for linear systems the error already at $X = 10^4$ does not exceed 0.01% and decays as $1/(X \ln X)$.

Perspectives for further research

1. Non-abelian case – proving asymptotics for non-abelian groups remains an open problem. A natural tool is Artin L-functions and GRH for them.

2. Systems with dependent Galois groups – classification of all possible values for subgroups of symmetric groups.

3. Deeper analysis of the role of GRH for non-abelian systems, as well as for abelian systems where GRH allows a transition from slow logarithmic convergence to fast power convergence.

Concluding remarks

The paper shows that the Hardy–Littlewood and Bateman–Horn singular series possess a stable structure, and the convergence speed of their partial products is determined solely by the Galois group and the value of the coefficient $C(\mathcal{F})$.

- Linear systems (trivial group) give superfast convergence $(O(1/(x\ln x)))$.

- Abelian nonlinear systems with $C = 0$ converge slowly without GRH $\left(o(1/\ln x)\right)$ and fast with GRH $(O(x^{-1/2+\varepsilon}))$.

- Systems with $C \neq 0$ are not admissible in the classical conjectures because their singular product diverges.

The results provide a rigorous theoretical foundation for computing singular series and refine the Bateman–Horn formula. They also open the way to further study of non-abelian systems and the role of the Generalized Riemann Hypothesis